\documentclass[11pt]{article}

\usepackage{amsmath,amsfonts,amssymb}
\usepackage{mathdots}
\usepackage{graphicx}
\usepackage[colorlinks,linkcolor=blue,citecolor=blue]{hyperref}
\usepackage{color}

\newtheorem{theorem}{Theorem}[section]

\newtheorem{definition}[theorem]{Definition}

\title{A note on the structure of expansive matrices in indefinite inner product spaces}

\author{A.C.M.~Ran\footnote{Department of Mathematics, Faculty of Science, VU Amsterdam, De Boelelaan
    1111, 1081 HV Amsterdam, The Netherlands
    and Research Focus: Pure and Applied Analytics, North-West~University,
Potchefstroom,
South Africa. E-mail:
    \texttt{a.c.m.ran@vu.nl}}}
    
\date{}

\begin{document}

    \maketitle

\noindent
{\bf Abstract}

\noindent
A result on the structure of expansive matrices in an indefinite inner product space is derived, which exhibits the largest unitary compression of the matrix.



\section{Introduction}

Matrices and operators with structure in an indefinite inner product space have received substantial attention in the last decades, see e.g., the books \cite{AI,B,IKL}. In the finite dimensional situation
much of that attention is focussed on deriving canonical forms, see e.g., \cite{GLR, GLR1} for $H$-selfadjoint matrices. For $H$-unitary matrices there is also a substantial literature on canonical forms, see e.g., \cite{GJR, FL} and the references cited therein. Several classes of normal matrices in indefinite inner product spaces are considered in \cite{Meh06a, Meh06}. For some classes of matrices canonical forms do not exist, but simple forms can be deduced. This is the case for $H$-dissipative matrices, see \cite{RT}, and for $H$-expansive matrices \cite{FGJR}.
In several cases canonical forms or even simple forms are either very hard or impossible, see \cite{MRR1,MRR2} for some results on $H$-hyponormal matrices. Most of the attention in the literature has been on forms that make it easy to establish the existence of invariant maximal $H$-semidefinite subspaces, even in the case where the matrix $H$ is not invertible and so the indefinite inner product is a degenerate one, see \cite{MRR}. 

As a starting point for this paper, let us mention some results from \cite{MRR}, Section 3, which we specify to the case where $H=H^*$ is invertible. A matrix $A$ is called $H$-expansive if $A^*HA-H=D\geq 0$. In that case it holds that $\ker A$ is $H$-negative (this result goes back to Ju.P. Ginsburg and is also stated in Chapter 11 of \cite{IKL}). In fact, more is true: there exists an invertible matrix $S$ such that 
$$
S^*HS=\begin{bmatrix} H_1 & 0 \\ 0 & -I\end{bmatrix}, S^*A^*HAS=\begin{bmatrix} C_1 &0 \\ 0 & 0\end{bmatrix},
$$
and $H_1$ and $C_1$ are both invertible, and $C_1-H_1 \geq 0$. This implies in particular that $\ker A^*HA$ is $H$-negative. This result may be viewed as a result on the structure of $H$-expansive matrices. In this note we shall establish a more detailed structure theorem for $H$-expansive matrices.


On $\mathbb{C}^n$ consider an indefinite inner product induced by the invertible Hermitian matrix $H$. Further, let $A$ be an $n\times n$ matrix, and let $N$ be an $A$-invariant subpace of $\mathbb{C}^n$. We start by outlining a construction that can be found at several places in the literature for several classes of matrices. Define  $M=N\cap (HN)^\perp$. Choose a complement $M_1$ of $M$ in $N$, then $M_1$ is $H$-nondegenerate. For each such $M_1$, there is an $M_2$ which is skewly-linked to $M$, $H$-neutral and $H$-orthogonal to $M_1$. Then with respect to the decompositon
$$\mathbb{C}^n=\underbrace{M\dot+ M_1}_{N} \dot+ M_2 \dot+ M_3$$ 
we have
$$
A=\begin{bmatrix} 
A_{11} & A_{12} & A_{13} & A_{14} \\ A_{21} & A_{22} & A_{23} & A_{24} \\ 0 & 0 & A_{33} & A_{34} \\ 0 & 0 & A_{43} & A_{44}
\end{bmatrix},
\qquad
H=\begin{bmatrix}
0 & 0 & I & 0 \\ 0 & H_{22} & 0 & 0 \\ I & 0 & 0 & 0 \\ 0 & 0  & 0 & H_{44}
\end{bmatrix}.
$$
In this case,
if $A$ and or $N$ have further properties with respect to the indefinite inner product defined by $H$, then there will be relations between the block entries $A_{ij}$ in the matrix $A$, and further properties of $H_{22}$ and $H_{44}$. 

For instance, take first the case where $A$ is $H$-selfadjoint, that is $HA=A^*H$. First, it will be shown that $M$ is $A$-invariant. Indeed, let $x\in M$, then $Ax\in N$ because $N$ is $A$-invariant. Further, for all $y\in  N$ we have $\langle HAx, y\rangle =\langle Hx, Ay\rangle=0$, since $x\in M=N\cap (HN)^\perp$ and $Ay\in N$. It follows that $Ax\in (HN)^\perp$, and so $Ax\in N\cap (HN)^\perp=M$. Hence in this case the $A$-invariance of $N$ also implies that $M$ is $A$-invariant, and so $A_{21}=0$. Writing out $HA=A^*H$ in terms of the blocks in the $4\times 4$ block matrices yields
$$
A=\begin{bmatrix} 
A_{11} & A_{12} & A_{13} & A_{14} \\ 0 & A_{22} & H_{22}^{-1}A_{12}^*& 0 \\ 0 & 0 & A_{11}^* & 0 \\ 0 & 0 & H_{44}^{-1}A_{14}^* & A_{44}
\end{bmatrix},
\qquad
H=\begin{bmatrix}
0 & 0 & I & 0 \\ 0 & H_{22} & 0 & 0 \\ I & 0 & 0 & 0 \\ 0 & 0  & 0 & H_{44}
\end{bmatrix},
$$
with $A_{22}$ being $H_{22}$-selfadjoint and $A_{44}$ being $H_{44}$-selfadjoint.\footnote{The author learned this decomposition of the space and the corresponding block forms of $A$ and $H$ for the $H$-selafdjoint case from H. Langer, private communication.} Clearly, in case $N$ is $H$-nonnegative, respectively, $H$-nonpositive, $H_{22}$ is positive definite, respectively, negative definite \cite{Langer}.

A similar construction, but in far less detail with respect to the submatrices $A_{ij}$, can be found in \cite{AI}, page 198, for the case of $H$-expansive matrices. In the next section we shall make this decomposition much more precise for $H$-expansive matrices.

The construction above is also reminiscent of the construction in the proofs of Theorem 6 and Proposition 4 in \cite{MRR2} for the case of $H$-hyponormal matrices. The special case where $N=\ker A$ and $A$ is $H$-normal is discussed in Lemma 2 in \cite{MRR3}.

\section{The structure theorem for $H$-expansive matrices}

Let $A$ be $H$-expansive, that is $A^*HA-H=D\geq 0$. Here $H=H^*$ is invertible, and $A$, $H$ and $D$ are $n\times n$ complex matrices. In this section we develop a structure result for $H$-expansive matrices.

\begin{definition}
Let $A$ be $H$-expansive. We shall say that a pair $(A_2,H_2)$ is a \emph{unitary compression} of $(A,H)$ if there is an $H$-nondegenerate subspace $K$ and an $H$-neutral subspace $M$ such that with respect to the decomposition $\mathbb{C}^n=M\dot+ K\dot+ (HK)^\perp$ we have
$$
A=\begin{bmatrix} A_1 & A_{12} & A_{13} \\ 0 & A_2 & A_{23} \\ 0 & 0 & A_3\end{bmatrix} , \qquad H=\begin{bmatrix} 0 & 0 & H_{13} \\ 0 & H_2 & 0 \\ H_{13}^* & 0 & H_{33}\end{bmatrix},
$$
and $A^*HA-H=\begin{bmatrix} 0 & 0 & 0 \\ 0 & 0 & 0 \\ 0 & 0 & \widetilde{D}\end{bmatrix}$ for some positive semidefinite $\widetilde{D}$.

If in addition the subspace $K$ is $A$-invariant (and so $A_{12}=0$) then we shall say that the pair $(A_2, H_2)$ is a \emph{unitary part} of $(A,H)$.
\end{definition}

The definition is partly motivated by a result in \cite{IKL}, to be precise, Lemma 11.5, which states for an $H$-expansive operator on a $\Pi_\kappa$ space, that if $L$ is an $A$-invariant subspace such that $AL=L$ and $\langle HAx,Ax\rangle = \langle Hx,x\rangle$ for all $x\in L$, then $(H^{-1}A^*H) L=L$, and $A|_L: L\to L$ is an $H$-unitary operator. It is also shown that if $L$ is an $H$-neutral $A$-invariant subspace then $AL=L$ and  $A|_L:L\to L$ is $H$-unitary.

Recall that the pair $(D,A)$ is called \emph{observable} if $\cap_{j=0}^\infty \ker DA^j=\{0\}$. In case $(D,A)$ is not observable, the subspace $\cap_{j=0}^\infty \ker DA^j$ is called the \emph{unobservable subspace} of the pair $(D,A)$.

The following statement is the main result of this paper. 

\begin{theorem}
Let $A$ be an $n\times n$ $H$-expansive matrix, with $A^*HA-H=D\geq 0$. Introduce the $(D,A)$-unobservable subspace $N=\cap_{j=0}^{n-1} \ker DA^j$, and set $M=N\cap (HN)^\perp$. Choose a complement $M_1$ of $M$ in $N$, then $M_1$ is $H$-nondegenerate. For each such $M_1$, there is an $M_2$ which is skewly-linked to $M$, $H$-neutral and $H$-orthogonal to $M_1$. Take $M_3=(H(M\dot+ M_1\dot+ M_2))^\perp$.
With respect to the decompositon $\mathbb{C}^n=M\dot+ M_1 \dot+ M_2 \dot+ M_3$ we have
$$
A=\begin{bmatrix} 
A_{11} & A_{12} & A_{13} & A_{14} \\ 0 & A_{22} & -H_{22}^{-1}A_{22}^{-*}A_{12}^*A_{11}^{-*} & 0 \\ 0 & 0 & A_{11}^{-*} & 0 \\ 0 & 0 & A_{43} & A_{44}
\end{bmatrix},
\qquad
H=\begin{bmatrix}
0 & 0 & I & 0 \\ 0 & H_{22} & 0 & 0 \\ I & 0 & 0 & 0 \\ 0 & 0  & 0 & H_{44}
\end{bmatrix},
$$
where $A_{22}$ is $H_{22}$-unitary and $A_{44}$ is $H_{44}$-expansive. Moreover, set
\begin{align*}
D_{11}
&= A_{11}^{-1}A_{13} +A_{13}^*A_{11}^{-*} +A_{11}^{-1}A_{12}H_{22}^{-1}A_{12}^*A_{11}^{-*}+A_{43}^*H_{44}A_{43}
\\
D_{12}&= A_{11}^{-1}A_{14} +A_{43}^*H_{44}A_{44} ,\\
D_{22}&=A_{44}^*H_{44}A_{44}-H_{44}.
\end{align*}
Then 
$$
A^*HA-H=\begin{bmatrix} 0 & 0 & 0 & 0 \\ 0 & 0 &0 & 0 \\ 0 & 0 & D_{11} & D_{12} \\ 0 & 0 & D_{12}^* & D_{22}\end{bmatrix} \geq 0.
$$
In addition, the pair 
$$
\left(\begin{bmatrix} D_{11} & D_{12} \\ D_{12}^* & D_{22}\end{bmatrix}, \begin{bmatrix} A_{11}^{-*} & 0 \\ A_{43} &A_{44}\end{bmatrix} \right)
$$ 
is observable. 
\end{theorem}

Note that the pair $(A_{22},H_{22})$ is a unitary compression of $(A,H)$ and by construction it is (up to similarity) the largest unitary compression of $(A,H)$. It is a unitary part of $(A,H)$ if $A_{12}=0$.

\bigskip

{\bf Proof.}
Note that $K$ and $M$ in the definition are contained in the kernel of $D$, and that their sum must be $A$-invariant. Hence it is natural to consider the largest $A$-invariant subpace of $\ker D$.
Thus, consider the subspace $N=\cap_{j=0}^{n-1}\ker( DA^j)$, so $N$ is the unobservable subspace for the pair $(D,A)$ (which is the largest $A$-invariant subspace in $\ker D$). Then $N$ is $A$-invariant, and for $x\in N$ we have $A^*HAx=Hx$. Since this is the case, $Ax\not=0$ for all $x\in N\setminus\{0\}$, and $A:N\to N$ is invertible. Hence also $AN=N$.

Let $M=N\cap (HN)^\perp$. We claim that $M$ is $A$-invariant. 
Indeed, suppose $x\in M$ and consider $Ax$. Since $x\in M$ we have $y^*Hx=0$ for all $y\in N$. Since $x\in N$ we have $Dx=0$, so $y^*Dx=0$ for all $y\in N$. It follows that $y^*A^*HAx=y^*(H+D)x=0$
for all $y\in N$. But since $AN=N$ we have that any $z\in N$ can be written as $z=Ay$ for some $y\in N$. Hence for all $z\in N$ we have $z^*HAx=y^*(H+D)x=0$. Thus $Ax$ is $H$-orthogonal to $N$, that is $Ax\in N\cap (HN)^\perp=M$. 

Choose a complement $M_1$ of $M$ in $N$. Then for all $x\in M_1$ and all $m\in M$ we have $\langle Hx,m\rangle =0$. Indeed, this is the case since $x\in M_1\subset N$ and $m\in M\subset (HN)^\perp$. 

Next, we show that $M_1$ must be $H$-nondegenerate. Suppose this is not the case, then there exists a vector $x\in M_1$ such that for all $y\in M_1$ we have $\langle Hx,y\rangle =0$. Take a vector $z\in N$ and write it as $z=m+y$, where $m\in M$ and $y\in M_1$. Then 
$\langle Hx,z\rangle = \langle Hx,m\rangle+\langle Hx,y\rangle =0$. But this means that $x\in N\cap (HN)^\perp=M$, which is a contradiction.

So $M_1$ is $H$-nondegenerate, $M_1\subset (HM)^\perp$ and for all $x\in M_1$ we have
$A^*HAx=Hx$. 

Take a subspace $M_2$ which is skewly-linked to $M$, and such that $M_2$ is $H$-neutral and $H$-orthogonal to $M_1$ (this can always be achieved). Then $M\dot+ M_1\dot+ M_2$ is $H$-nondegenerate. Let $M_3=(H(M\dot+ M_1\dot+ M_2))^\perp$. Then 
$$
\mathbb{C}^n= \underbrace{M\dot+ M_1}_{N} \dot+ M_2 \dot+ M_3.
$$
Using that $N$ and $M$ are both $A$-invariant, and using the properties of the subspaces with respect to the indefinite inner product,
we can write $A$ and the Gramm matrix for $H$ as follows
$$
A=\begin{bmatrix} 
A_{11} & A_{12} & A_{13} & A_{14} \\ 0 & A_{22} & A_{23} & A_{24} \\ 0 & 0 & A_{33} & A_{34} \\ 0 & 0 & A_{43} & A_{44}
\end{bmatrix},
\qquad
H=\begin{bmatrix}
0 & 0 & I & 0 \\ 0 & H_{22} & 0 & 0 \\ I & 0 & 0 & 0 \\ 0 & 0  & 0 & H_{44}
\end{bmatrix}.
$$
We shall consider the submatrices in the third and fourth column of $A$ to develop some structure. 

For all $y$ and for all $x\in N$ we write these vectors as
$$
y=\begin{bmatrix} y_1 \\ y_2 \\ y_3\\\ y_4 \end{bmatrix} , \qquad x=\begin{bmatrix} x_1 \\ x_2 \\ 0 \\ 0 \end{bmatrix} .
$$
Since $A^*HAx=Hx$ for $x\in N$, we have $\langle HAx, Ay\rangle =\langle Hx,y\rangle $. Writing this in terms of the coordinates we have
$$
\left\langle H\begin{bmatrix} A_{11}x_1+A_{12}x_2 \\ A_{22}x_2 \\ 0 \\ 0 \end{bmatrix}, Ay\right\rangle =
\left\langle \begin{bmatrix} 0 \\ H_{22} A_{22}x_2 \\  A_{11}x_1+A_{12}x_2 \\ 0 \end{bmatrix}, Ay\right\rangle=
\langle H_{22} x_2, y_2\rangle + \langle x_1, y_3\rangle.
$$
This can be worked out a bit further to
\begin{align*}
&\langle H_{22} A_{22}x_2 ,A_{22}y_2 \rangle + 
\langle H_{22} A_{22}x_2 ,A_{23}y_3+A_{24}y_4 \rangle +
\langle A_{11}x_1+A_{12}x_2 , A_{33}y_3+A_{34}y_4\rangle\\ & =
\langle H_{22} x_2, y_2\rangle + \langle x_1, y_3\rangle.
\end{align*}
Using that $x_2 , y_2\in N$ we have $\langle H_{22} A_{22}x_2 ,A_{22}y_2 \rangle =\langle H_{22} x_2, y_2\rangle $, i.e., $A_{22}$ is $H_{22}$-unitary. That means we have
$$
\langle H_{22} A_{22}x_2 ,A_{23}y_3+A_{24}y_4 \rangle +
\langle A_{11}x_1+A_{12}x_2 , A_{33}y_3+A_{34}y_4\rangle =
 \langle x_1, y_3\rangle.
$$
Now take $x_2=0$ and $y_3=0$ to see that for all $x_1\in M$ and all $y_4\in M_3$ we have
$\langle A_{11}x_1 , A_{34}y_4\rangle =0$. Since $M$ and $M_3$ are skewly linked, and $A_{11}$ is invertible (remember that $A|_N:N\to N$ is invertible) this implies that $A_{34}=0$. Next, take $x_2=0$ and $y_4=0$, to obtain $\langle A_{11}x_1 , A_{33}y_3\rangle =
 \langle x_1, y_3\rangle$. Recall that ${\rm dim\,}M={\rm dim\,}M_3$. Choosing appropriate bi-orthogonal bases in $M$ and $M_3$ we can identify operators on these spaces as matrices of the same size. Then the equality $\langle A_{11}x_1 , A_{33}y_3\rangle =
 \langle x_1, y_3\rangle$ implies that $A_{33}^*A_{11}=I$. Then we have that
$$
\langle H_{22} A_{22}x_2 ,A_{23}y_3+A_{24}y_4 \rangle +
\langle A_{12}x_2 , A_{33}y_3\rangle =0.
$$
Now take $y_3=0$ to obtain $A_{24}=0$, and then the remaining equation is
$$
\langle H_{22} A_{22}x_2 ,A_{23}y_3 \rangle +
\langle A_{12}x_2 , A_{33}y_3\rangle =0.
$$
Since this must hold for all choices of $x_2$ and $y_3$ we have
$$
A_{23}^*H_{22}A_{22} =-A_{33}^*A_{12}=-A_{11}^{-1}A_{12}.
$$
This determines $A_{23}$. 

We arrive at the following structure result: with respect to the decomposition
$$
\mathbb{C}^n=M\dot+ M_1 \dot+ M_2 \dot+ M_3
$$
we have that $A$ and $H$ are given by
$$
A=\begin{bmatrix} 
A_{11} & A_{12} & A_{13} & A_{14} \\ 0 & A_{22} & -H_{22}^{-1}A_{22}^{-*}A_{12}^*A_{11}^{-*} & 0 \\ 0 & 0 & A_{11}^{-*} & 0 \\ 0 & 0 & A_{43} & A_{44}
\end{bmatrix},
\qquad
H=\begin{bmatrix}
0 & 0 & I & 0 \\ 0 & H_{22} & 0 & 0 \\ I & 0 & 0 & 0 \\ 0 & 0  & 0 & H_{44}
\end{bmatrix},
$$
with $A_{22}$ being $H_{22}$-unitary. 

There are further restrictions, which we obtain by using that $A^*HA-H=D$ should be positive semidefinite.

Computing $A^*HA$ we obtain, using $A_{22}^*H_{22}A_{22}=H_{22}$ and the formula for $A_{23}$
\begin{align*}
&\begin{bmatrix}
A_{11}^* & 0 & 0 & 0 \\
A_{12}^* & A_{22}^* & 0 & 0 \\ A_{13}^* & A_{23}^* & A_{11}^{-1} & A_{43}^* \\
A_{14}^* & 0 & 0 & A_{44}^*
\end{bmatrix}
\begin{bmatrix}
0 & 0 & I & 0 \\ 0 & H_{22} & 0 & 0 \\ I & 0 & 0 & 0 \\ 0 & 0  & 0 & H_{44}
\end{bmatrix}
\begin{bmatrix} 
A_{11} & A_{12} & A_{13} & A_{14} \\ 0 & A_{22} & A_{23} & 0 \\ 0 & 0 & A_{11}^{-*} & 0 \\ 0 & 0 & A_{43} & A_{44}
\end{bmatrix}\\
=&
\begin{bmatrix} 
0 & 0 & I & 0 \\ 0 & H_{22} & 0 & 0 \\ I & 0 & D_{11} & D_{12} \\
0 & 0 & D_{21} & D_{22}+H_{44}
\end{bmatrix},
\end{align*}
where
\begin{align*}
D_{11}&= A_{11}^{-1}A_{13} +A_{13}^*A_{11}^{-*} +A_{11}^{-1}A_{12}A_{22}^{-1}H_{22}^{-1}A_{22}^{-*}A_{12}^* A_{11}^{-*}+A_{43}^*H_{44}A_{43}\\
&= A_{11}^{-1}A_{13} +A_{13}^*A_{11}^{-*} +A_{11}^{-1}A_{12}H_{22}^{-1}A_{12}^*A_{11}^{-*}+A_{43}^*H_{44}A_{43}
\\
D_{12}&= A_{11}^{-1}A_{14} +A_{43}^*H_{44}A_{44} ,\\
D_{21}&=D_{12}^*,\\
D_{22}&=A_{44}^*H_{44}A_{44}-H_{44}.
\end{align*}
So, 
$$
A^*HA-H=\begin{bmatrix} 0 & 0 & 0 & 0 \\ 0 & 0 &0 & 0 \\ 0 & 0 & D_{11} & D_{12} \\ 0 & 0 & D_{12}^* & D_{22}\end{bmatrix}
$$
with $D_{11}, D_{12}$ and $D_{22}$ as above. Since $A^*HA-H$ is positive semidefinite, we have in particular that $A_{44}$ is $H_{44}$-expansive.

Moreover, from the construction of $N$ we have that the pair
$$
\left(\begin{bmatrix} D_{11} & D_{12} \\ D_{12}^* & D_{22}\end{bmatrix}, \begin{bmatrix} A_{11}^{-*} & 0 \\ A_{43} &A_{44}\end{bmatrix} \right)
$$ 
is observable. 
\hfill$\Box$

\bigskip

The first three examples below will show that the dimension of $M_1$ strongly depends on the rank of $A^*HA-H$.

{\bf Example 1.} Consider 
$$
A=J_5(1)=\begin{bmatrix} 1 & 1 & 0 & 0 & 0 \\ 0 & 1 & 1 & 0 & 0 \\ 0 & 0 & 1 & 1 & 0 \\ 0 & 0 & 0 & 1 & 1 \\ 0 & 0 & 0 & 0 & 1 \end{bmatrix}, \qquad
H=\begin{bmatrix} 0 & 0 & 0 & 0 & 1 \\ 0 & 0 & 0 & -1 & -1 \\ 0 & 0 & 1 & 2 & 2 \\ 0 & -1 & 2 & 4 & 6 \\ 1 & -1 & 2 & 6 & 2\end{bmatrix}.
$$
One ckecks that $A$ is $H$-expansive, in fact
$$
D=A^*HA-H= \begin{bmatrix} 0 & 0 & 0 & 0 & 0 \\  0 & 0 & 0 & 0 & 0 \\  0 & 0 & 0 & 0 & 0 \\ 
0 & 0 & 0 & 5 & 8 \\ 0 & 0 & 0 & 8 & 16\end{bmatrix} \geq 0.
$$
Then the $(D,A)$-unobservable subspace $N$ is given by $N={\rm span\,} \{e_1, e_2, e_3\}$, and $M$ is given by
$M={\rm span\,}\{e_1, e_2\}$. Set $M_1={\rm span\,}\{e_3\}$, and
$$
x_1=\begin{bmatrix} 3 & -2 & 0 & -1 & 1 \end{bmatrix}^T, \quad 
x_2=\begin{bmatrix} 0 & 0 & 2 & -1 & 0\end{bmatrix}.
$$
Then $M_2={\rm span\, } \{x_1,x_2\}$ is skewly-linked to $M$, $H$-orthogonal to $M_1$ and $H$-neutral. Then $\mathbb{C}^5=M\dot+ M_1\dot+ M_2$. Take as basis $\{e_1, e_2, e_3, x_1, x_2\}$, and let $S$ be the matrix with these vectors as its columns (in the given order). We have
$$
S^{-1}AS=\begin{bmatrix} 1 & 1 & 0 & -2 & 0 \\ 0 & 1 & 1 & 0 & 2 \\ 0 & 0 & 1 & 1 & -1 \\
0 & 0 & 0 & 1 & 0 \\ 0 & 0 & 0 & -1 & 1 \end{bmatrix}, \quad
S^*HS=
\begin{bmatrix} 0 & 0 & 0 & 1 & 0 \\ 0 & 0 & 0 & 0 & 1 \\ 0 & 0 & 1 & 0 & 0 \\
1 & 0 & 0 & 0 & -4 \\ 0 & 1 & 0 & -4 & 0
\end{bmatrix}.
$$
So the largest unitary compression is formed by the pair $A_{22}=1, H_{22}=1$ appearing in the center of the matrices $S^{-1}AS$ and $S^*HS$, respectively. \hfill$\Box$

{\bf Example 2.}
Consider the same matrix $A=J_5(1)$ as in the previous example, but now with the indefinite inner product given by 
$$
H=\begin{bmatrix} 0 & 0 & 0 & 0 & 1 \\ 0 & 0 & 0 & -1 & \tfrac{3}{2} \\ 0 & 0 & 1 & -\tfrac{1}{2} & \tfrac{1}{2} \\ 0 & -1 &  -\tfrac{1}{2} & 0 & 1 \\ 1 &  \tfrac{3}{2}  &  \tfrac{1}{2}  & 1 & 0\end{bmatrix}.
$$
Then $A^*HA-H={\rm diag\, }(0, 0, 0, 0, 2)$, so $A$ is $H$-expansive. Observe that $N=\cap_{j=0}^4 \ker DA^j= {\rm span\, }\{ e_1, e_2, e_3, e_4\}$ and $M={\rm span\, }\{e_1\}$. Set $M_1={\rm span\, }\{ e_2, e_3, e_4\}$.  To determine $M_2$, which has to be one-dimensional, write $M_2={\rm span\, }\{x\}$, and solve for the coordinates of $x$ from the conditions that $M_2$ should satisfy. We obtain
\begin{align*}
&\langle He_1,x\rangle =1 \mbox{ giving } x_5=1,\\
&\langle He_2,x\rangle =0 \mbox{ giving } x_4=\tfrac{3}{2},\\
&\langle He_3,x\rangle =0 \mbox{ giving } x_3=\tfrac{1}{4},\\
&\langle He_4,x\rangle =0 \mbox{ giving } x_2=\tfrac{7}{8},\\
&\langle H x,x\rangle =0 \mbox{ giving } x_1=-\tfrac{47}{32}.\\
\end{align*}
Set $S=\begin{bmatrix} e_1 & e_2 & e_3 & e_4 & x\end{bmatrix}$, then 
$$
S^{-1}AS=\begin{bmatrix}1 & 1 & 0 & 0 & \tfrac{7}{8}\\
0 & 1 & 1 & 0 & \tfrac{1}{4}\\
0 & 0 & 1 & 1 & \tfrac{3}{2}\\
0 & 0 & 0 & 1 & 1 \\ 0 & 0 & 0 & 0 & 1\end{bmatrix}, \qquad
S^*HS=\begin{bmatrix} 0 & 0 & 0 & 0 & 1 \\ 0 & 0 &0  & -1 & 0 \\ 0 & 0 & 1 & -\tfrac{1}{2} & 0 \\ 
0 & -1 & -\tfrac{1}{2} & 0 & 0 \\
1 & 0 & 0 & 0 & 0\end{bmatrix}.
$$
It follows that the largest unitary compression is the pair
$$
A_{22}=\begin{bmatrix}1 & 1 & 0 \\
0 & 1 & 1 \\0 & 0 & 1\end{bmatrix}, \qquad 
H_{22}=\begin{bmatrix} 0 & 0 & -1 \\ 0 & 1 & -\tfrac{1}{2} \\ -1 & -\tfrac{1}{2} & 0
\end{bmatrix}.
$$
 \hfill$\Box$
 
{\bf Example 3.} Next, we consider the same matrix $A=J_5(1)$, but now with the indefinite inner product given by 
$$
H=\begin{bmatrix} 0 & 0 & 0 & 0 & 1 \\ 0 & 0 & 0 & -1 & \tfrac{3}{2} \\ 0 & 0 & 1 & -\tfrac{1}{2} & \tfrac{1}{2} \\ 0 & -1 &  -\tfrac{1}{2} & 0 & 0 \\ 1 &  \tfrac{3}{2}  &  \tfrac{1}{2}  & 0 & 0\end{bmatrix}.
$$
Then $A$ is $H$-unitary, in fact the pair $(A,H)$ is in the canonical form given in \cite{GJR}. We have $N=\ker D=\mathbb{C}^5$, $M=(0)$ and hence $M_1$ is the whole space. As expected, the largest unitary compression in the sense of our definition is the pair $(A,H)$ itself, and it is clearly also the largest unitary part in this case. \hfill$\Box$

The next example is an example of a pair $(A,H)$ with a non-trivial unitary part.

{\bf Example 4.}
Let $A=J_2(1)\oplus J_2(1)\oplus J_2(1)$
and
$$
H=\begin{bmatrix} 0 & 0 & 0 & 1 \\
0 & 0 & -1 & 0 \\
0 & -1 & 0 & 0 \\
1 & 0 & 0 & 0 \end{bmatrix}
\oplus
\begin{bmatrix} 0 & 1  \\
1 & 0\end{bmatrix} .
$$
Then $A$ is $H$-expanisve, and $D=A^*HA-H= 0_{5\times 5}\oplus (2)$. So $N={\rm span\, }
\{e_1, e_2, e_3, e_4, e_5\}$ and $M={\rm span\,} \{e_5\}$. Take $M_1={\rm span\,}\{e_1, e_2, e_3, e_4\}$. Clearly, $M_2$ can be taken to be $M_2={\rm span\, }\{e_6\}$. So in this case the largest unitary compression if the pair formed by the $4\times 4$ matrices
$$
A_{22}=J_2(1)\oplus J_2(1), H_{22}=\begin{bmatrix} 0 & 0 & 0 & 1 \\
0 & 0 & -1 & 0 \\
0 & -1 & 0 & 0 \\
1 & 0 & 0 & 0 \end{bmatrix}.
$$
In fact, this is the largest unitary part in this case. \hfill$\Box$

The final example is one without a unitary compression.

{\bf Example 5.} Let $A=J_4(1)$ and 
$$
H=\begin{bmatrix} 0 & 0 & 0 & 1 \\ 0 & 0 & -1 & -1\\ 0 & -1 & 0 & 0 \\ 1 & -1 & 0 & 0\end{bmatrix}.
$$
Then $D=A^*HA-H={\rm diag\,}(0,0,2,0)$ and so $N={\rm span\,}\{e_1, e_2\}$. Since this space is a $H$-neutral, we have $M=N$, and hence $M_1=\{0\}$.

\paragraph{Acknowledgements:}
This work is based on research supported in part by the National Research Foundation of South Africa (Grant Number 145688).

The author wishes to thank Christian Mehl for valuable discussions and suggestions.


\begin{thebibliography}{w}
\bibitem{AI}
{ T. Ya. Azizov, I.S. Iohvidov},
{\it Linear operators in spaces
with an indefinite metric}, John Wiley and Sons, Chicester, 1989
[Russian original 1979].


\bibitem{B}
{ J. Bogn\'ar},
{\it Indefinite inner product spaces.} Ergebnisse der Mathematik und ihrer Grenzgebiete, Band 78, Springer-Verlag, New York-Heidelberg, 1974.

\bibitem{FL}
{ A. Ferrante, B.C. Levy}. {\rm Canonical form of symplectic matrix pencils,}
{\it Linear Algebra and its Applications}, 274:259-300, 1998.


\bibitem{FGJR}
J.H. Fourie, G.J. Groenewald, D.B. Janse van Rensburg, and A.C.M. Ran:
Simple forms and invariant subspaces of $H$-expansive matrices.
{\it Linear Algebra Appl.} 470 (2015), 300--340.

\bibitem{GLR}
{ I. Gohberg, P. Lancaster, L. Rodman},
{\it Matrices and Indefinite Scalar Products},
Oper. Theory: Adv. and Appl. 8, Birkh\"auser Verlag, Basel, 1983.

\bibitem{GLR1}
{I.~Gohberg, P.~Lancaster and L.~Rodman}.
{\it Indefinite Linear Algebra and Applications}, 
Birkh\"auser Verlag, Basel, 2005.

\bibitem{GJR}
{ G.J. Groenewald, D.B. Janse van Rensburg, A.C.M. Ran.}
{\rm A canonical form for $H$-unitary matrices,} {\it Operators and Matrices.} 10(4): 739-783, 2016.

\bibitem{IKL}
{ I.S. Iohvidov, M.G. Krein, H. Langer},
{\it Introduction to the spectral theory of operators in spaces with an indefinite metric},
Mathematical Research 9, Akademie-Verlag, Berlin, 1982.

\bibitem{Langer}
H. Langer. Private communication. 
\bibitem{Meh06a}
{ Chr. Mehl. }
{\rm On classification of normal matrices in indefinite inner product spaces,}
{\it Electronic Journal of Linear Algebra}, 15:50--83, 2006.

\bibitem{Meh06}
{ Chr. Mehl.}
{\rm Essential decomposition of polynomially normal matrices in real indefinite inner product spaces,}
{\it Electronic Journal of Linear Algebra}, 15:84--106, 2006.

\bibitem{MRR}
Chr. Mehl, A.C.M. Ran, L. Rodman. Semidefinite invariant subspaces: degenerate inner products.
{\it Proceedings IWOTA 2002} 
{\it Oper. Theory Adv. Appl.} 149 (2004), 475--494.

\bibitem{MRR1}
Chr. Mehl, A.C.M. Ran and L. Rodman:
Hyponormal matrices and semidefinite invariant subspaces in
indefinite inner products.
{\it Electronic Journal of Linear Algebra} 11 (2004), 192-204.

\bibitem{MRR3}
Chr. Mehl, A.C.M. Ran and L. Rodman:
Polar decomposition of normal operators in indefinite inner
product spaces.
{\it Proceedings of 3d Workshop on Indefinite Inner Products},
{\it Oper. Theory Adv. Appl.}162, (2006) 277-292.

\bibitem{MRR2}
Chr. Mehl, A.C.M. Ran and L. Rodman:
Extension to maximal semidefinite invariant subspaces for
hyponormal matrices in indefinite inner products.
{\it Linear Algebra Appl.} 421 (2007), 110-116.

\bibitem{RT}
A.C.M. Ran and D. Temme:
Invariant semidefinite subspaces of dissipative matrices in an
indefinite inner product space, existence, construction and uniqueness.
{\it Linear Algebra Appl.} 212/213(1994), 169-214.
\end{thebibliography}
\end{document}